\documentclass[12pt]{article}
\usepackage{amsthm}
\usepackage{amsfonts}
\usepackage{epsfig}
\usepackage{amssymb}
\usepackage{amsmath}
\usepackage{amscd}
\usepackage{colonequals}
\usepackage{url}
\usepackage{float}
\newtheorem{theorem}{Theorem}
\newtheorem{lemma}[theorem]{Lemma}
\newtheorem{observation}[theorem]{Observation}

\newcommand{\BB}{\mathcal{B}}
\newcommand{\brm}[1]{\operatorname{#1}}

\title{Treewidth of graphs with balanced separations}

\author{Zden\v{e}k Dvo\v{r}\'{a}k\thanks{Charles University, Prague, Czech Republic.
E-mail: {\tt rakdver@iuuk.mff.cuni.cz}.  Supported by the Center of Excellence -- Inst. for Theor. Comp. Sci., Prague (project P202/12/G061 of Czech Science Foundation).}\and
Sergey Norin\thanks{Department of Mathematics and Statistics, McGill University. Email: {\tt snorin@math.mcgill.ca}}}
\date{}
\begin{document}
\maketitle

\begin{abstract}
We prove that if every subgraph of a graph $G$ has a balanced separation of order at most $a$ then $G$ has treewidth at most $15a$. This establishes a linear dependence between the treewidth and the separation number.
\end{abstract}

\section{Introduction and results}
Treewidth of a graph is an important graph parameter introduced by Robertson
and Seymour~\cite{rs2}, which expresses how ``tree-like" a given graph is. The
relation between the treewidth and other graph parameters, e.g. the maximum
order of a tangle~\cite{rs10} and the size of a largest
grid-minor~\cite{RSey,quickly,ChChGrid}, has been explored in a number of
papers; see~\cite{HaWoTree} for a recent survey. The goal of this paper is to
establish a linear dependence between treewidth and another parameter, the
separation number, which we now define.

A \emph{separation} of a graph $G$ is a pair $(A,B)$ of subsets of $V(G)$ such
that $A \cup B = V(G)$ and no edge of $G$ has one end in $A\setminus B$ and the other in
$B\setminus A$. The \emph{order} of the separation $(A,B)$ is  $|A \cap B|$. 
For an assignment $w:V(G)\to\mathbb{R}^+_0$ of non-negative weights to vertices of $G$
and a set $X\subseteq V(G)$, we define $w(X)=\sum_{v\in X} w(v)$.
We say that a separation $(A,B)$ of a graph $G$ is \emph{$w$-balanced} if $w(A\setminus B)\le \tfrac{2}{3}w(V(G))$
and $w(B\setminus A)\le \tfrac{2}{3}w(V(G))$.
A separation $(A,B)$ of a graph $G$ on $n$ vertices is \emph{balanced} if
$|A\setminus B|\le 2n/3$ and $|B\setminus A|\le 2n/3$, i.e., it is $w$-balanced for the function $w$
assigning weight~$1$ to each vertex of $G$.  The \emph{separation
number} $\brm{sn}(G)$ of a graph $G$ is the smallest integer $s$ such that every
subgraph of $G$ has a balanced separation of order at most $s$.

Let $\brm{tw}(G)$ denote the treewidth\footnote{We will not need the
definition of treewidth in this paper.  Instead, we use an equivalent parameter, the maximum order
of a bramble, which is defined in Section~\ref{sec-aux}.}
of the graph $G$.  The relation between the separation number and the treewidth has been
explored starting with Robertson and Seymour~\cite{rs2}.  They proved that
if $G$ has treewidth at most $t$, then for every assignment $w$ of non-negative weights
to vertices of $G$, there exists a $w$-balanced separation of $G$ of order at most $t+1$.
Conversely, they also showed that if $G$ has a $w$-balanced separation of order at most $s$
for every $\{0,1\}$-valued weight assignment $w$, then the treewidth of $G$ is less than $4s$.

Since $\brm{tw}(G')\le \brm{tw}(G)$ for every subgraph $G'$ of $G$, it follows that
$\brm{sn}(G) \leq \brm{tw}(G)+1$ for every graph $G$.  In this paper, we study whether
a weak converse to this inequality holds.  Note that this is not obvious: the fact that
each subgraph has a balanced separation of small order does not straightforwardly imply
that the whole graph has a $w$-balanced separation of small order for every
$\{0,1\}$-valued weight assignment $w$.

It is easy to show that $\brm{tw}(G) \leq 1 + \brm{sn}(G)\log(|V(G)|)$, see e.g.~\cite{bodlaender1998}.
Fox~\cite{Fox} stated without
proof that $\brm{sn}(G)$ and $\brm{tw}(G)$ are tied, that is $\brm{tw}(G) \leq
f(\brm{sn}(G))$ for some function $f$. Finally, B{\"o}ttcher et
al.~\cite{BPTWBand} investigated the relation between the separation
number and the treewidth for graphs with bounded maximum degree. They have
shown that for fixed $\Delta$ and a subgraph-closed class of graphs with maximum
degree at most $\Delta$, the treewidth is sublinear in the number of vertices
if and only if the separation number is.

The following is our main result.

\begin{theorem}\label{thm-main}
The treewidth of any graph $G$ is at most $15\brm{sn}(G)$.
\end{theorem}

Note that Theorem~\ref{thm-main} implies the aforementioned result of
B{\"o}ttcher et al. without the restriction on the maximum degree.
We give the
necessary definitions and prove the preliminary results in
Section~\ref{sec-aux}. The proof of Theorem~\ref{thm-main} is completed
in Section~\ref{sec-proof}. 

\section{Brambles and flows}\label{sec-aux}

In this section we introduce the tools and auxiliary results necessary for the proof of Theorem~\ref{thm-main}.

\subsection{Brambles}

Let $G$ be a graph.
A set $\BB$ of non-empty subsets of $V(G)$ is a \emph{bramble} if the induced subgraph $G[S]$ is connected for all $S\in \BB$
and the induced subgraph of $G[S_1\cup S_2]$ is connected for all $S_1,S_2\in\BB$.
A \emph{hitting set} for a bramble $\BB$ is a set that intersects all elements of $\BB$.
The \emph{order} of the bramble is the smallest size of its hitting set.
The relation between the brambles and the treewidth of a graph is captured in the following result of Seymour and Thomas.

\begin{theorem}[Seymour and Thomas~\cite{bramble}]\label{thm-bramble}
A graph $G$ contains a bramble of order $k$ if and only if $G$ has treewidth at least $k-1$.
\end{theorem}

A \emph{flap assignment} of order $k$ in a graph $G$ is a function $\beta$ that maps each set $W\subseteq V(G)$ of size less than $k$
to the vertex set of a (non-empty) connected component of $G-W$.
A \emph{haven} of order $k$ is a flap assignment of order $k$ such that for all sets $W,W'\subseteq V(G)$ of sizes less than $k$,
the induced subgraph $G[\beta(W)\cup\beta(W')]$ is connected.  Clearly, if $\beta$ is a haven of order $k$ in $G$,
then $\{\beta(W):W\subseteq V(G),|W|<k\}$ forms a bramble in $G$ of order at least $k$.

Conversely, consider a bramble $\BB$ in $G$ of order $k$.
If a set $W\subseteq V(G)$ has size smaller than $k$, then $W$ is not a hitting set for $\BB$,
and thus there exists a set $S\in\BB$ disjoint from $W$.  Since $G[S]$ is connected, there exists a component $C$ of $G-W$ containing $S$.
Furthermore, since $G[S_1\cup S_2]$ is connected for all $S_1,S_2\in\BB$,
all elements of $\BB$ disjoint from $W$ are subsets of this unique component $C$.  Hence, setting $\beta(W)$ to be the vertex set
of such component $C$ for each $W\subseteq V(G)$ of size less than $k$ gives a haven $\beta$ of order $k$, which we call \emph{the haven
defined by the bramble $\BB$}.

Let us note a basic property of havens.
\begin{observation}\label{obs-haven}
Let $\beta$ be a haven of order $k$ in a graph $G$, let $W_1$ and $W_2$ be sets of vertices of $G$ of size less than $k$, and
let $(A,B)$ be a separation of $G$ of order less than $k$.
If $W_1\subseteq W_2$, then $\beta(W_2)\subseteq \beta(W_1)$.  Also, if\/ $\beta(W_1)\subseteq A\setminus B$,
then $\beta(A\cap B)\subseteq A\setminus B$.
\end{observation}

We need the following result enabling us to find a small(er) set in the image of a haven.

\begin{lemma}\label{lemma-shrink}
Let $a, r\ge 0$ be integers and let $G$ be a graph with separation number at most $a$.
Let $Y$ be a subset of vertices of $G$, and let $\beta$ be a haven in $G$ of order greater than $|Y|+ar$.
Then there exists a set $Y'$ of vertices of $G$ such that $|Y'|\le |Y|+ar$ and $|\beta(Y')|\le (2/3)^r|\beta(Y)|$.
\end{lemma}
\begin{proof}
We prove the claim by induction on $r$; for $r=0$, the claim clearly holds, hence assume that $r\ge 1$.
Let $G_1$ be the subgraph of $G$ induced by $\beta(Y)$.  Since $\brm{sn}(G)\le a$, there exists a balanced separation $(A,B)$ of $G_1$ of order
at most $a$.  Let $Y_1=Y\cup (A\cap B)$.  By Observation~\ref{obs-haven}, we have $\beta(Y_1)\subseteq \beta(Y)=V(G_1)$.
Since $\beta(Y_1)$ is disjoint from $A\cap B$ and induces a connected subgraph, we can by symmetry assume that
$\beta(Y_1)\subseteq A\setminus B$.  Since the separation $(A,B)$ is
balanced, we have $|\beta(Y_1)|\le |A\setminus B|\le \frac{2}{3}|\beta(Y)|$.  Since $|Y_1|\le |Y|+a$, the claim of the lemma
follows by the induction hypothesis applied to $Y_1$ and $r-1$.
\end{proof}

\subsection{Flows and clouds}
Let $\vec{G}$ be a directed graph and let $g: V(\vec{G}) \to \mathbb{R}_0^+$ be a function specifying the amount of flow entering its vertices. The value $g(V(\vec{G}))\colonequals\sum_{v\in V(\vec{G})}g(v)$ is called
the \emph{total supply}. 
Let $W=\{v_1,v_2,\ldots,v_k\} \subseteq V(\vec{G})$ be a set of sinks in $\vec{G}$, where a \emph{sink} is a vertex with no outgoing edges. A \emph{flow} towards $W$ in $\vec{G}$
with supply $g$ is a function $f:E(\vec{G}) \to \mathbb{R}_0^{+}$ satisfying the flow conservation equation
$$ g(v) + \sum_{e = uv \in E(\vec{G})}f(e) = \sum_{e = vz \in E(\vec{G})}f(e)$$
for every $v \in V(\vec{G}) \setminus W$. For a fixed flow $f$ and a vertex $v \in V(\vec{G})$, define $\brm{in}(v)$ to be the sum of the flow values on the edges into $v$.
The \emph{congestion at $v$} is defined as $\brm{in}(v)+g(v)$. The \emph{congestion of a flow $f$} is the maximum congestion over all the vertices of the graph.
Note that the total supply equals the sum of congestions of vertices of $W$.  Let $\Vert f\Vert=\sum_{e\in E(\vec{G})} f(e)$ and let the support $\brm{supp}(f)$ of $f$ be the set
of edges of $\vec{G}$ with non-zero value of the flow.

Let $G$ be an undirected graph and let $W$ be a subset of $V(G)$.  The graph $\vec{G}_W$ is obtained from $G$ by replacing each edge with a pair of oppositely
directed edges and deleting the edges directed away from the vertices of $W$.
A \emph{$W$-cloud} in $G$ is a flow in $\vec{G}_W$ towards $W$ such that the supply at each vertex is at most $1$.
The following result is an easy corollary of the max-flow min-cut theorem (applied to the graph obtained from $\vec{G}_W$ by splitting each vertex
into an input and an output part connected by an edge of capacity $\alpha$, adding a source vertex with edges of capacity $1$ going into all
input parts of vertices of $\vec{G}_W$, and having the capacity of all the other edges infinite).

\begin{observation}\label{obs-cloud}
Let $G$ be an undirected graph and let $W$ be a subset of $V(G)$, and let $\alpha,s>0$ be real numbers.
There exists a $W$-cloud in $G$ with congestion at most $\alpha$ and total supply greater than $s$ if and only if
$G$ contains no separation $(C,D)$ with $W\subseteq C$ and $|C\setminus D|+\alpha|C\cap D|\le s$.
\end{observation}

Let $G$ be an undirected graph and let $W$ be a subset of $V(G)$.  For a $W$-cloud $f$ with supply $g$, a vertex $v\in V(G)$ is \emph{saturated}
if $g(v)=1$, and $v$ is \emph{hungry} if $v$ is not saturated, but has non-zero congestion or belongs to $W$.  Let $\brm{sat}(f)$ denote the set of saturated vertices of $G$.
Let $\vec{G}_W\downarrow f$ denote the subgraph of $\vec{G}_W$ whose vertices are all hungry and saturated vertices of $\vec{G}$ and with edge set equal to the support of $f$,
and let $G_W\downarrow f$ denote the underlying undirected graph of $\vec{G}_W\downarrow f$.  We say that a $W$-cloud $f$ is \emph{$\alpha$-tame}
if the congestion of $f$ is at most $\alpha$ and additionally,
\begin{itemize}
\item[(i)] the hungry vertices have no in-neighbors in $\vec{G}_W\downarrow f$,
\item[(ii)] each saturated vertex has at most one hungry in-neighbor in $\vec{G}_W\downarrow f$, and
\item[(iii)] if $\alpha\le 1$, then $\vec{G}_W\downarrow f$ is the edgeless graph with vertex set $W$.
\end{itemize}

\begin{lemma}\label{lemma-mincloud}
Let $G$ be an undirected graph, let $W$ be a subset of $V(G)$, and let $\alpha,s>0$ be real numbers.
If there exists any $W$-cloud in $G$ with congestion at most $\alpha$ and total supply exactly $s$,
then there also exists an $\alpha$-tame $W$-cloud with total supply exactly $s$.
\end{lemma}
\begin{proof}
Let $f$ be a $W$-cloud in $G$ with congestion at most $\alpha$ and total supply exactly $s$ such that
among all such $W$-clouds, $\Vert f\Vert$ is minimum, and subject to that, $|\brm{supp}(f)|-|\brm{sat}(f)|$ is minimum.
Note that such a $W$-cloud exists by the compactness of the space of all $W$-clouds with given total supply and upper bound on congestion.

If $\alpha\le 1$, then let $f'$ be the flow with the empty support and with supply at each vertex of $W$ equal to the congestion of $f$ at $W$.
Since $\Vert f\Vert\le \Vert f'\Vert=0$, we have $f=f'$ and $\vec{G}_W\downarrow f$ is the edgeless graph with vertex set $W$.  Clearly,
$f$ is $\alpha$-tame.

Hence, assume that $\alpha>1$, and in particular the condition (iii) of $\alpha$-tameness is trivially satisfied. The graph $\vec{G}_W\downarrow f$ is acyclic,
since if $C$ were a directed cycle in $\vec{G}_W\downarrow f$, we could decrease the values of the flow on the edges of $C$ by the same positive amount and obtain a $W$-cloud $f'$
with congestion at most $\alpha$ and total supply $s$ such that $\Vert f'\Vert<\Vert f\Vert$.

Let $g$ be the supply of $f$.  Suppose first that there exists an edge $uv\in E(\vec{G}_W\downarrow f)$ such that $v$ is hungry.  Let $P$ be a maximal directed path in $\vec{G}_W\downarrow f$
ending with the edge $uv$.  Let $\delta_1>0$ be the minimum of the flow values on the edges of $P$.
Since $\vec{G}_W\downarrow f$ is acyclic, the starting vertex $x$ of $P$ has no incoming edges, and thus $g(x)\ge \delta_1$.  Let $\delta=\min(\delta_1, 1-g(v))$.
Let $f'$ be obtained from $f$ by decreasing the flow by $\delta$ on the edges of $P$, decreasing the supply at $x$ by $\delta$ and increasing the supply at $v$ by $\delta$.
Clearly, $f'$ is a $W$-cloud with the same total supply and at most as large congestion as $f$, and $\Vert f'\Vert<\Vert f\Vert$, which is a contradiction.
Hence, hungry vertices have no in-neighbors in $\vec{G}_W\downarrow f$, and $f$ satisfies the condition (i) of $\alpha$-tameness.

Suppose now that there exist distinct edges $v_1z,v_2z\in E(\vec{G}_W\downarrow f)$ such that $v_1$ and $v_2$ are hungry.
By the previous paragraph, $v_1$ has no incoming edges, and thus $g(v_1)\ge f(v_1z)$.
Let $\delta=\min(f(v_1z),1-g(v_2))$.
Let $f'$ be the flow obtained by decreasing the flow on $v_1z$ by $\delta$, increasing the flow on $v_2z$ by $\delta$,
decreasing the supply at $v_1$ by $\delta$ and increasing the supply at $v_2$ by $\delta$.
Only the congestion of $v_2$ increases, and since $v_2$ has no in-neighbors by the previous paragraph,
the congestion of $f'$ at $v_2$ is at most its supply at $v_2$, which is at most $1\le \alpha$.
Hence, $f'$ is a $W$-cloud with the same total supply as $f$ and congestion at most $\alpha$, and $\Vert f'\Vert=\Vert f\Vert$.
Furthermore, either $\brm{supp}(f')=\brm{supp}(f)\setminus\{v_1z\}$, or $\brm{sat}(f')=\brm{sat}(f)\cup \{v_2\}$, and
thus $|\brm{supp}(f')|-|\brm{sat}(f')|<|\brm{supp}(f)|-|\brm{sat}(f)|$, which is a contradiction.
We conclude that each vertex has at most one hungry in-neighbor in $\vec{G}_W\downarrow f$, and $f$ satisfies the condition (ii) of $\alpha$-tameness.

Therefore, $f$ is $\alpha$-tame.
\end{proof}

We now use clouds to extract a special subgraph of a hypothetical counterexample to Theorem~\ref{thm-main}.
\begin{lemma}\label{lemma-nicesg}
Let $G$ be a graph and let $a=\brm{sn}(G)$.  If $a\ge 1$ and $G$ has treewidth at least $15a$,
then $G$ contains a subgraph $F$ with a separation $(X,Z)$ of order $14a$ such that the following conditions hold.
Let $W=X\cap Z$, $s=14|Z\setminus W|$ and $\alpha=\tfrac{s}{7a}$.
\begin{itemize}
\item[\textnormal{(a)}] The subgraph $F[X]$ contains an $\alpha$-tame $W$-cloud $f$ with total supply greater than $s$,
such that each vertex of $F[X]$ is hungry or saturated.
\item[\textnormal{(b)}] There exists a balanced separation $(A,B)$ of $F$ such that
$$|A\cap (W\cup B)\cap X|\le a.$$
\end{itemize}
\end{lemma}
\begin{proof}
By Theorem~\ref{thm-bramble}, $G$ contains a bramble $\BB$ of order at least $15a+1$.
Let $\beta$ denote the haven of order at least $15a+1$ defined by $\BB$.
Let $W$ be a subset of vertices of $G$ of size at most $14a$ chosen so that $|\beta(W)|$ is minimum, and subject to that $|W|$ is minimum.
Clearly, $|W|=14a$, as otherwise for any $w\in \beta(W)$, we have $\beta(W\cup \{w\})\subsetneq \beta(W)$, contradicting the minimality of $|\beta(W)|$.
Similarly, the minimality of $|W|$ implies that every vertex of $W$
has a neighbor in $\beta(W)$.  Let $X'=V(G)\setminus \beta(W)$ and $Z=W\cup\beta(W)$.
Let $s=14|Z\setminus W|=14|\beta(W)|$ and $\alpha=\tfrac{s}{7a}$.

Suppose first that $G[X']$ has a separation $(C,D)$ such that $W\subseteq C$ and $|C\setminus D|+\alpha|C\cap D|\le s$,
and thus the order of $(C,D)$ is at most $s/\alpha=7a$ and $|C\setminus D|\le s=14|\beta(W)|$.
Let $C'=C\cup \beta(W)$ and $Y=C'\cap D=C\cap D$, and note that $(C',D)$ is a separation of $G$.
Since $\beta(W)\subseteq C'\setminus D$, Observation~\ref{obs-haven} implies 
$\beta(Y)\subseteq C'\setminus D$.  Consequently,
$|\beta(Y)|\le |C'\setminus D|=|\beta(W)|+|C\setminus D|\le 15|\beta(W)|$.
Applying Lemma~\ref{lemma-shrink} with $r=7$, we conclude that there exists a set $Y'\subseteq V(G)$
of size at most $|Y|+ra\le 14a$ such that $|\beta(Y')|\le (2/3)^7|\beta(Y)|\le 15(2/3)^7|\beta(W)|<|\beta(W)|$,
which contradicts the choice of $W$.

Therefore, $G[X']$ has no such separation, and thus by Observation~\ref{obs-cloud}, $G[X']$ contains
a $W$-cloud $f$ of congestion at most $\alpha$ and total supply $s_f>s$.
Furthermore, by Lemma~\ref{lemma-mincloud}, we can assume that $f$ is $\alpha$-tame.
Let $X$ be the set of hungry and saturated vertices of $G[X']$ (so $X=V(G[X']_W\downarrow f)$)
and let $F=(G[X']_W\downarrow f)\cup G[Z]$.  Clearly, $(X,Z)$ is a separation of $F$, $X\cap Z=W$ has size $14a$,
and $f$ is an $\alpha$-tame $W$-cloud in $F[X]$ satisfying (a).

Since $\brm{sn}(G)=a$ and $F$ is a subgraph of $G$, there exists a balanced separation $(A,B)$ of $F$ of order at most $a$.
Let $(A',B')=(A\cap Z,B\cap Z)$ be the separation of $G[Z]$ corresponding to $(A,B)$.
If $B'\setminus A'$ contains all vertices of $\beta(W)$, then $W\subseteq B'$, since each vertex of $W$ has a neighbor in $\beta(W)$.
Consequently $W\subseteq B$ and $A\cap (W\cup B)\cap X\subseteq A\cap B$, implying (b).
Hence, we can assume that $\beta(W)\not\subseteq B'\setminus A'$, and symmetrically (swapping the labels of $A$ and $B$),
we have $\beta(W)\not\subseteq A'\setminus B'$.

Let $W'=W\cup (A'\cap B')$. 
We have $|W'|\le 15a$, and thus $\beta(W')$ is defined.  Since $W\subseteq W'$, we have $\beta(W')\subseteq \beta(W)$,
and by symmetry (swapping the labels of $A$ and $B$ if necessary), we can assume that
$\beta(W')\subseteq B'\setminus (A'\cup W)\subsetneq \beta(W)$.
Let $W''=W'\cap B'$.  No vertex of $W\setminus B'$ has
a neighbor in $B'\setminus A'$, and thus $\beta(W'')=\beta(W')$.  Consequently, the minimality of $|\beta(W)|$
implies that $|W''|>14a=|W|$, and thus $|A\cap (W\cup B)\cap X|=|A\cap B|+|A\cap W|-|A\cap B\cap Z|=|A\cap B|+|W|-|W''|<|A\cap B|\le a$.
Therefore, (b) holds.
\end{proof}

We are now ready to prove Theorem~\ref{thm-main} by showing that the conditions (a) and (b) are contradictory.

\section{Proof of Theorem~\ref{thm-main}}\label{sec-proof}
Let $a=\brm{sn}(G)$. We assume without loss of generality that $a \geq 1$.
Suppose for a contradiction that $G$ has treewidth at least $15a$.  By Lemma~\ref{lemma-nicesg},
$G$ contains a subgraph $F$ with a separation $(X,Z)$ of order $14a$ satisfying
(with $W=X\cap Z$, $s=14|Z\setminus W|$ and $\alpha=\tfrac{s}{7a}$) conditions (a) and (b).
Let $(A,B)$ be a balanced separation of $F$ such that
\begin{equation}\label{eq-smside}
|A\cap (W\cup B)\cap X|\le a,
\end{equation}
which exists by (b).

By (a), $F[X]$ contains an $\alpha$-tame $W$-cloud $f$ with total supply greater than $s$,
such that each vertex of $F[X]$ is hungry or saturated, and thus $X=V(F[X]_W\downarrow f)$.
If $\alpha<1/2$, then $X=W$ by the condition (iii) of $\alpha$-tameness.
Consequently, $V(F)=Z$.  Furthermore, $|Z\setminus W|=s/14=\alpha a/2<a/4$. 
By (\ref{eq-smside}), we have $|A\cap W|\le a$, and thus
$$
\frac{|B\setminus A|}{|V(F)|}\ge \frac{|W|-|A\cap W|}{|W|+|Z\setminus W|}>\frac{13a}{(14+1/4)a}>2/3.$$
This contradicts the fact that $(A,B)$ is a balanced separation of $F$, and thus $\alpha\ge 1/2$.

Let $g$ be the supply of the $W$-cloud $f$.  Let $a_0=|A\cap (W\cup B)\cap X|$ and let $a_1\le a_0$ be the number of hungry vertices
in $A\cap (W\cup B)\cap X$.  Note that by the condition (i) of $\alpha$-tameness, no flow enters the hungry vertices of $A\cap (W\cup B)\cap X$
through incident edges.  Consider a saturated vertex $v\in A\cap X$.  The flow $g(v)=1$ originating at $v$ must reach $W$, and
thus it must pass through (or end in) a saturated vertex of $A\cap (W\cup B)\cap X$.  Since all vertices have congestion at most $\alpha$,
it follows that the number of saturated vertices of $A\cap X$ is at most $\alpha (a_0-a_1)$.  By the condition (ii) of $\alpha$-tameness,
each saturated vertex of $A\cap X$ has at most
one hungry in-neighbor, and each hungry vertex $x$ of $A\cap X$ either belongs to $A\cap (W\cup B)\cap X$ or has a saturated out-neighbor in $A\cap X$.
Since $a_0\le a$ by (\ref{eq-smside}), we conclude that
$$|A|\le |A\cap X|+|Z\setminus W|\le 2\alpha (a_0-a_1)+a_1+|Z\setminus W|\le 2\alpha a+|Z\setminus W|=\frac{5}{14}s.$$
On the other hand, the total supply of $f$ is greater than $s$ and each vertex has supply at most $1$,
giving the lower bound greater than $s$ on the size of $X$.  Hence, we have
$$|V(F)|=|X|+|Z\setminus W|>s+|Z\setminus W|=\frac{15}{14}s.$$
Consequently,
$$
\frac{|B\setminus A|}{|V(F)|}=1-\frac{|A|}{|V(F)|}>1-\frac{5s/14}{15s/14}=2/3.
$$
This contradiction to the fact that $(A,B)$ is a balanced separation of $F$ finishes the proof. 

\section*{Acknowledgments}

We would like to thank the anonymous referee whose suggestions helped us to significantly improve the multiplicative constant in Theorem~\ref{thm-main}.
An earlier version of the paper~\cite{septw-v1} used results of Chen, Kleinberg, Lov\'{a}sz, Rajaraman, Sundaram, and Vetta~\cite{confluent} regarding confluent flows; while we were
able to avoid their usage eventually, they inspired the method we use.


\end{document}